\documentclass[12pt,leqno,draft]{article}

\newtheorem{theorem}{Theorem}
\newtheorem{lemma}[theorem]{Lemma}
\newtheorem{proposition}[theorem]{Proposition}
\newtheorem{definition}[theorem]{Definition}
\newtheorem{corollary}[theorem]{Corollary}

\newcommand{\begintheorem}{\addtocounter{equation}{1}\begin{theorem}}
\newcommand{\beginlemma}{\addtocounter{equation}{1}\begin{lemma}}
\newcommand{\beginproposition}{\addtocounter{equation}{1}\begin{proposition}}
\newcommand{\begindefinition}{\addtocounter{equation}{1}\begin{definition}}
\newcommand{\begincorollary}{\addtocounter{equation}{1}\begin{corollary}}

\begin{document}

\title{Some topics in complex and harmonic analysis}

\author{Stephen William Semmes	\\
	Rice University		\\
	Houston, Texas}

\date{}

\maketitle

	Fix a positive integer $n$, and let $C_b({\bf R}^n)$
denote the vector space of continuous complex-valued functions
on ${\bf R}^n$ which are bounded.  If $f \in C_b({\bf R}^n)$,
then the supremum norm of $f$ is given by
\begin{equation}
	\|f\| = \sup \{|f(x)| : x \in {\bf R}^n\}.
\end{equation}
Observe that
\begin{equation}
	\|f_1 + f_2\| \le \|f_1\| + \|f_2\|
\end{equation}
and
\begin{equation}
	\|f_1 \, f_2\| \le \|f_1\| \, \|f_2\|
\end{equation}
for all $f_1, f_2 \in C_b({\bf R}^n)$.

	By a finite measure on ${\bf R}^n$ we mean a linear mapping
$\lambda$ from $C_b({\bf R}^n)$ into the complex numbers such that
there is a nonnegative real number $L$ with the property that
\begin{equation}
\label{|lambda(f)| le L ||f||}
	|\lambda(f)| \le L \, \|f\|
\end{equation}
for all $f \in C_b({\bf R}^n)$ and $\lambda$ is continuous in the
sense that if $\{f_j\}_{j=1}^\infty$ is a sequence of bounded
continuous functions on ${\bf R}^n$ which is uniformly bounded and
converges uniformly on compact subsets of ${\bf R}^n$ to a bounded
continuous function $f$, then $\{\lambda(f_j)\}_{j=1}^\infty$
converges to $\lambda(f)$.  If $f$ is any bounded continuous function
on ${\bf R}^n$, then there is a sequence $\{f_j\}_{j=1}^\infty$ of
continuous functions on ${\bf R}^n$ such that each $f_j$ has compact
support in ${\bf R}^n$, the $f_j$'s are uniformly bounded, and the
$f_j$'s converge to $f$ uniformly on compact subsets of ${\bf R}^n$.
As a result, a finite measure on ${\bf R}^n$ is determined by its
restriction to the vector space of continuous functions with compact
support on ${\bf R}^n$.  In fact, if one starts with a linear
functional on the vector space of continuous functions on ${\bf R}^n$
with compact which is bounded in the sense of the same inequality
(\ref{|lambda(f)| le L ||f||}), then it admits a unique extension
to a finite measure on ${\bf R}^n$.

	If $\lambda$ is a finite measure on ${\bf R}^n$, then its norm
$\|\lambda\|_*$ is defined to be the supremum of $|\lambda(f)|$ over
all $f \in C_b({\bf R}^n)$ such that $\|f\| \le 1$.  This is
equivalent to defining $\|\lambda\|_*$ to be the supremum of
$|\lambda(f)|$ over all continuous functions $f$ on ${\bf R}^n$ with
compact support such that $\|f\| \le 1$.  It is also the same as the
smallest nonnegative real number $L$ such that (\ref{|lambda(f)| le L
||f||}) holds.  Notice that 
\begin{equation}
	\|\lambda_1 + \lambda_2\|_* \le \|\lambda_1\|_* + \|\lambda_2\|_*
\end{equation}
for all finite measures $\lambda_1$, $\lambda_2$ on ${\bf R}^n$.

	As a basic class of examples, suppose that $h(x)$ is a
continuous function on ${\bf R}^n$ which is integrable in the sense
that $\int_{{\bf R}^n} |h(x)| \, dx$ is finite.  This leads to a
finite measure $\lambda$ on ${\bf R}^n$ defined by setting
$\lambda(f)$ equal to $\int_{{\bf R}^n} h(x) \, f(x) \, dx$ for all
bounded continuous functions $f$ on ${\bf R}^n$.  To be a bit more
precise, these integrals can be defined as improper integrals, which
reduce to classical Riemann integrals when $f(x)$ has compact support.
The norm of this linear functional is equal to $\int_{{\bf R}^n}
|h(x)| \, dx$.

	Now suppose that $A$ is a subset of ${\bf R}^n$ which is at
most countable, and that $a(x)$ is a complex-valued function on $A$
such that $\sum_{x \in A} |a(x)|$ is finite.  This leads to a finite
measure $\lambda$ on ${\bf R}^n$ defined by setting $\lambda(f)$
to be equal to $\sum_{x \in A} a(x) \, f(x)$.  The norm of this
linear functional is equal to $\sum_{x \in A} |a(x)|$.  One can also
define finite measures on ${\bf R}^n$ by integrating over submanifolds
of ${\bf R}^n$ of any dimension.

	If $\lambda$ is a finite measure on ${\bf R}^n$ and $\phi$ is
a bounded continuous function on ${\bf R}^n$, then we can get a new
finite measure $\lambda_\phi$ on ${\bf R}^n$ by putting
$\lambda_\phi(f) = \lambda (\phi \, f)$ for all bounded continuous
functions $f$ on ${\bf R}^n$.  It is easy to check that
$\|\lambda_\phi\|_* \le \|\phi\| \, \|\lambda\|_*$.  One can show that
if $\{\phi_j\}_{j=1}^\infty$ is a sequence of bounded continuous
functions on ${\bf R}^n$ which are uniformly bounded and converge to
$1$ uniformly on compact subsets of ${\bf R}^n$, then the
corresponding $\lambda_{\phi_j}$'s converge to $\lambda$ in norm,
which is to say that $\|\lambda_{\phi_j} - \lambda\|_* \to 0$ as $j
\to \infty$.  In particular, one can choose the $\phi_j$'s so that
they have compact support in ${\bf R}^n$.  In other words, finite
measures on ${\bf R}^n$ can be approximated in norm by measures with
compact support.

	Suppose that $\mu$, $\nu$ are finite measures on ${\bf R}^m$,
${\bf R}^n$, respectively.  We can define a new finite measure $\mu
\times \nu$ on ${\bf R}^m \times {\bf R}^n$, which can be identified
with ${\bf R}^{m + n}$, by saying that the value of $\mu \times \nu$
applied to a bounded continuous function $F(x, y)$ on ${\bf R}^m \times
{\bf R}^n$ is obtained first by integrating $F$ in $x$ using $\mu$
to get a bounded continuous function of $y$, and then integrating
that using $\nu$ to get a complex number.  Of course one could also
apply $\mu$ and $\nu$ in the other order.  Either way one gets the
same result, because the two approaches give the same answer
when $F(x, y)$ is of the form $f_1(x) \, f_2(y)$ where $f_1$, $f_2$
are bounded continuous functions on ${\bf R}^m$, ${\bf R}^n$,
which is to say that they both are equal to $\mu(f_1) \, \nu(f_2)$.
By linearity they both give the same answer when $F(x, y)$ is a 
finite linear combinations of products of functions of $x$, $y$
separately, and one can conclude that they give the same answer
for all bounded continuous functions $F(x, y)$ through suitable
approximation arguments.

	To be a bit more precise there are some subtleties here
concerning the fact that $\mu \times \nu$ is continuous with respect
to uniformly bounded sequences of functions which converge uniformly
on compact subsets.  One might prefer to start by defining $\mu \times
\nu$ applied to continuous functions $F(x, y)$ which have compact
support.  At any rate, because the finite measures $\mu$, $\nu$ can be
approximated in norm by measures with compact support, it is easy to
see that everything works fine.  Notice that $\mu \times \nu$ has
compact support if $\mu$, $\nu$ do.  One can also check that the norm
of $\mu \times \nu$ is equal to the product of the norms of $\mu$,
$\nu$ for any finite measures $\mu$, $\nu$.

	If $\lambda$ is a finite measure on ${\bf R}^n$ and $f$ is a
bounded continuous function on ${\bf R}^n$, then we define the
convolution $\lambda * f$ to be the function on ${\bf R}^n$ obtained
by setting $(\lambda * f)(x)$ equal to $\lambda(f_x)$, where $f_x$ is
the bounded continuous function on ${\bf R}^n$ given by $f_x(y) = f(x
- y)$.  One can check that $\lambda * f$ is continuous, and it is
bounded with supremum norm less than or equal to the product of the
norm of $\lambda$ and the supremum norm of $f$.  If $\mu$, $\nu$ are
two finite measures on ${\bf R}^n$, then we define their convolution
$\mu * \nu$ to be the finite measure obtained by setting $(\mu *
\nu)(f)$ for a bounded continuous function $f$ on ${\bf R}^n$ equal to
$(\mu \times \nu)(F)$, $F(x, y) = f(x + y)$.  The norm of $\mu * \nu$
is less than or equal to the product of the norms of $\mu$ and $\nu$.

	Let $\lambda$ be a finite measure on ${\bf R}^n$.  For each
$\xi \in {\bf R}^n$, put $e_\xi(x) = \exp ( - 2 \pi i \, \xi \cdot x)$,
where $\xi \cdot x = \sum_{j=1}^n \xi_j \, x_j$ is the usual inner
product of vectors in ${\bf R}^n$.  This is a bounded continuous
function on ${\bf R}^n$, with $|e_\xi(x)| = 1$ for all $x, \xi \in
{\bf R}^n$.  Define the Fourier transform of $\lambda$ by
$\widehat{\lambda}(\xi) = \lambda(e_\xi)$ for all $\xi \in {\bf R}^n$.
One can check that this is a continuous function on ${\bf R}^n$, and
in fact that it is uniformly continuous.  Also, $\widehat{\lambda}$ is
bounded, with supremum norm less than or equal to the norm of
$\lambda$.

	Suppose that $\mu$, $\nu$ are finite measures on ${\bf R}^n$,
so that their Fourier transforms are bounded continuous functions on
${\bf R}^n$.  The \emph{multiplication formula} states that
$\mu(\widehat{\nu})$ is equal to $\nu(\widehat{\mu})$.  Indeed, each
of these is equal to $\mu \times \nu$ applied to the function $F(x,
\xi) = \exp (- 2 \pi i \, \xi \cdot x)$.  Notice too that the Fourier
transform of the convolution $\mu * \nu$ is equal to the product of
the Fourier transforms of $\mu$, $\nu$.  If $\lambda$ is a finite
measure on ${\bf R}^n$, then $\lambda * e_\xi(x)$ is equal to
$\widehat{\lambda}(-\xi)$ times $e_\xi(x)$.

	If $z$, $\zeta$ are elements of ${\bf C}^n$, which is to say
that they are $n$-tuples of complex numbers, then we can define $z
\cdot \zeta$ in the same manner as before, as $\sum_{j=1}^n z_j \,
\zeta_j$.  We can define $e_\zeta(z)$ for $z, \zeta \in {\bf C}^n$
through the same formula as before, i.e., $e_\zeta(z) = \exp (- 2 \pi
i \, \zeta \cdot z)$.  This is a complex analytic function of
$z$ and $\zeta$.  For instance, one can expand it out into a power
series in the $z_j$'s and $\zeta_j$'s, using the usual power series
expansion for the exponential.

	If $A$ is a nonempty closed subset of ${\bf R}^n$, let us
define $\widehat{A}$ to be the subset of ${\bf C}^n$ consisting of the
$\zeta \in {\bf C}^n$ such that $e_\zeta(x)$ is bounded on $A$.  For
$\zeta \in \widehat{A}$, let us put $a(\zeta)$ equal to the supremum
of $|e_\zeta(x)|$ over $x \in A$.  Thus ${\bf R}^n \subseteq
\widehat{A}$ and $a(\xi) = 1$ for all $\xi \in {\bf R}^n$.  In
particular, $\widehat{A}$ is not empty.  If $A$ is bounded, then
$\widehat{A} = {\bf C}^n$.

	In general, if $\zeta = \xi + i \eta$ with $\xi, \eta \in {\bf
R}^n$, then $\zeta \in \widehat{A}$ if and only if $\eta \in {\bf
R}^n$, and in this case $a(\zeta) = a(i \eta)$.  If $\zeta_1, \zeta_2
\in \widehat{A}$, then $\zeta_1 + \zeta_2 \in \widehat{A}$, and
$a(\zeta_1 + \zeta_2)$ is less than or equal to the product of
$a(\zeta_1)$ and $a(\zeta_2)$.  If $\zeta \in \widehat{A}$ and $t$ is
a nonnegative real number, then $t \, \zeta \in \widehat{A}$ and $a(t
\, \zeta)$ is equal to $a(\zeta)^t$.

	It follows that $\widehat{A}$ is actually a tube over a convex
cone.  To be more precise, let $A^*$ denote the set of $\eta \in {\bf
R}^n$ such that $i \eta \in \widehat{A}$.  Clearly $0 \in A^*$,
$\eta_1 + \eta_2 \in A^*$ when $\eta_1, \eta_2 \in A^*$, and $t \,
\eta \in A^*$ when $\eta \in A^*$ and $t$ is a nonnegative real
number.  In other words, $A^*$ is a convex cone, and $\widehat{A}$ is
equal to the set of $\zeta \in {\bf C}^n$ of the form $\zeta = \xi + i
\eta$ with $\xi, \eta \in {\bf R}^n$ and $\eta \in A^*$.

	Now suppose that $\lambda$ is a finite measure on ${\bf R}^n$
with support contained in $A$, in the sense that $\lambda(f) = 0$
whenever $f$ is a bounded continuous function on ${\bf R}^n$ such that
$f(x) = 0$ for all $x \in A$.  In this event we can define
$\lambda(f)$ for all bounded continuous functions on $A$, by extending
any such function to a bounded continuous function on ${\bf R}^n$ and
applying $\lambda$ to the extension.  The value of $\lambda$ applied
to the extension of $f$ does not depend on the choice of the
extension, because $\lambda$ applied to a function that vanishes on
$A$ is $0$.  Note that $|\lambda(f)|$ is less than or equal to the
norm of $\lambda$ times the supremum of $|f(x)|$, $x \in A$, for any
bounded continuous function $f$ on $A$, because a bounded continuous
extension of $f$ to all of ${\bf R}^n$ can always be chosen so that
the supremum norm of the extension is less than or equal to the
supremum of $|f(x)|$, $x \in A$.

	For each $\zeta \in \widehat{A}$, $e_\zeta(x)$ defines a
bounded continuous function on $A$, and thus we can extend the Fourier
transform to $\widehat{A}$ by putting $\widehat{\lambda}(\zeta) =
\lambda(e_\zeta)$.  Thus we get that $|\widehat{\lambda}(\zeta)|$ is
less than or equal to the norm of $\lambda$ times $a(\zeta)$ for all
$\zeta \in \widehat{A}$.  For $\zeta \in {\bf C}^n$ such that $-\zeta
\in \widehat{A}$ and $x \in {\bf R}^n$ we have that $e_\zeta(x - y)$
is a bounded continuous function of $y$ on $A$, and one can define
$(\lambda * e_\zeta)(x)$ to be equal to $\lambda$ applied to $e_\zeta
(x - y)$ as a function of $y$ as before.  Once again we also have that
$(\lambda * e_\zeta)(x)$ is equal to $\widehat{\lambda}(-\zeta)$
times $e_\zeta(x)$.

	Consider for the moment the special case where $A$ is bounded.
If $\lambda$ is a finite measure on ${\bf R}^n$ with support contained
in $A$, then we can define $\lambda(f)$ in a natural way for any
continuous function $f$ on ${\bf R}^n$.  Namely, if $f$ is not
bounded, we can replace it with a bounded continuous function which is
equal to it on $A$, and for that matter we can replace it with a
continuous function on ${\bf R}^n$ with compact support which is equal
to $f$ on $A$.  Of course $f$ is automatically bounded on $A$, since
$A$ is closed and bounded and therefore compact.

	In this situation $\widehat{A} = {\bf C}^n$ and
$\widehat{\lambda}(\zeta) = \lambda(e_\zeta)$ is defined for all
$\zeta \in {\bf C}^n$.  The Fourier transform of $\lambda$ is in fact
a complex analytic function on ${\bf C}^n$.  One way to look at this
is that the Fourier transform of $\lambda$ is a smooth function on
${\bf C}^n$, and its restriction to any complex line in ${\bf C}^n$ is
a complex analytic function of a single complex variable.  One can
also look at this in terms of a power series expansion for the Fourier
transform of $\lambda$ which converges on all of ${\bf C}^n$.  At any
rate, this complex analytic extension of the Fourier transform to
${\bf C}^n$ is uniquely determined by its restriction to ${\bf R}^n$.

	For each $\zeta \in {\bf C}^n$ we have that
$|\widehat{\lambda}(z)| \le \|\lambda\|_* \, a(\zeta)$.  More
precisely, if $\zeta = \xi + i \eta$, with $\xi, \eta \in {\bf R}^n$,
then $a(\zeta) = a(i \eta)$, and thus $|\widehat{\lambda}(\zeta)| \le
\|lambda\|_* \, a(i \eta)$.  In particular $\widehat{\lambda}(\xi + i
\eta)$ is bounded as a function of $\xi$ for each fixed $\eta$.  We
can also describe $a(i \eta)$ as the maximum of $\exp (2 \pi \eta
\cdot x)$ over $x \in A$.  Notice especially that $a(i \eta)$
is bounded by the exponential of a constant times the norm of $\eta$.

	If $\alpha$ is a positive real number, define the functions
$G_\alpha(z)$ and $W_\alpha(z)$ on ${\bf C}^n$ by $\exp (- 4 \pi^2
\alpha \, z \cdot z)$ and $(4 \pi \alpha)^{-n/2} \exp (- z \cdot z /
(4 \alpha))$, respectively.  These are complex analytic functions on
${\bf C}^n$ whose restrictions to ${\bf R}^n$ are the usual
Gauss--Weierstrass kernels.  On ${\bf R}^n$ these functions are
integrable, and hence have associated finite measures with these
functions as densities.  The Fourier transforms of these functions are
defined to be the Fourier transforms of the associated measures, and
it is well known that the Fourier transforms of these functions are
equal to each other.  Actually, these functions have sufficient decay
on ${\bf R}^n$ so that the Fourier transforms may be defined on all of
${\bf C}^n$, and the Fourier transforms of $G_\alpha$ and $W_\alpha$
are equal to each other on all of ${\bf C}^n$.

	Let $\lambda$ be a finite measure on ${\bf R}^n$, and for each
$\alpha > 0$ and $x \in {\bf R}^n$ put $\psi_{\alpha, x}(\xi) =
G_\alpha(\xi) \exp(2 \pi i \, x \cdot \xi)$.  Thus $\psi_{\alpha,
x}(xi)$ is an integrable function of $\xi$, so that its Fourier
transform is defined.  By the multiplication formula, $\int_{{\bf
R}^n} \widehat{\lambda}(\xi) \, \psi_{\alpha, x}(\xi) \, d\xi$ is
equal to $\lambda$ applied to the Fourier transform of $\psi_{\alpha,
x}$.  Because the Fourier transform of $G_\alpha$ is equal to
$W_\alpha$, the Fourier transform of $\psi_{\alpha, x}$ evaluated at
some point $y \in {\bf R}^n$ is equal to $W_\alpha(y - x)$.

	In other words, $\int_{{\bf R}^n} \widehat{\lambda}(\xi) \,
G_\alpha(\xi) \, \exp(2 \pi i \, x \cdot \xi) \, d\xi$ is equal to
$(\lambda * W_\alpha)(x)$.  For each $\alpha > 0$ we can identify
$W_\alpha$ on ${\bf R}^n$ with a finite measure, since $W_\alpha$ is
integrable, and thus for each bounded continuous function $f$ on ${\bf
R}^n$ we can define $W_\alpha * f$ as a bounded continuous function on
${\bf R}^n$.  The convolution $\lambda * W_\alpha$ is actually
integrable on ${\bf R}^n$, and for each bounded continuous function
$f$ on ${\bf R}^n$ we have that $\int_{{\bf R}^n} (\lambda *
W_\alpha)(x) \, f(x) \, dx$ is equal to $\lambda$ applied to $W_\alpha
* f$.  This can be shown through standard arguments.

	If $f$ is a bounded continuous function on ${\bf R}^n$, then
the convolutions $W_\alpha * f$ are uniformly bounded and converge to
$f$ as $\alpha \to 0$ uniformly on compact subsets of ${\bf R}^n$.
One might say that $\lambda * W_\alpha$ tends to $\lambda$ as $\alpha
\to 0$ in a weak sense, which is that the integral of $\lambda *
W_\alpha$ times a bounded continuous function $f$ tends to
$\lambda(f)$ as $\alpha \to 0$.  Indeed, these integrals are equal
to $\lambda$ applied to $W_\alpha * f$, and the latter are
uniformly bounded and converge to $f$ uniformly on compact subsets
of ${\bf R}^n$.

	As a result we obtain that a finite measure on ${\bf R}^n$ is
uniquely determined by its Fourier transform.  Namely, if $\lambda_1$,
$\lambda_2$ are two finite measures on ${\bf R}^n$ such that
$\widehat{\lambda}_1(\xi) = \widehat{\lambda_2}(\xi)$ for all $\xi \in
{\bf R}^n$, then $\lambda = \lambda_1 - \lambda_2$ is a finite measure
on ${\bf R}^n$ such that $\widehat{\lambda}(\xi) = 0$ for all $\xi \in
{\bf R}^n$.  It follows from the preceding discussion that $(\lambda *
W_\alpha)(x) = 0$ for all $x \in {\bf R}^n$, and hence that
$\lambda(f) = 0$ for all bounded continuous functions $f$ on ${\bf
R}^n$.

	We have seen that if $\lambda$ is a finite measure on ${\bf
R}^n$ which is supported in a compact set, then the Fourier transform
of $\lambda$ extends to a complex analytic function on all of ${\bf
C}^n$.  If we also assume that the Fourier transform of $\lambda$
vanishes on a nonempty open subset of ${\bf R}^n$, then it follows
from well known results in complex analysis that the Fourier transform
of $\lambda$ vanishes everywhere.  In particular, a finite measure
$\lambda$, on ${\bf R}^n$ with compact support whose Fourier transform
has compact support as a function on ${\bf R}^n$ is equal to $0$,
i.e., $\lambda(f) = 0$ for all bounded continuous functions $f$ on
${\bf R}^n$.

	Let $\lambda$ be a finite measure on ${\bf R}^n$, so that the
Fourier transform $\widehat{\lambda}(\xi)$ defines a bounded
continuous function on ${\bf R}^n$.  Suppose that the Fourier
transform of $\lambda$ is integrable, which is to say that $\int_{{\bf
R}^n} |\widehat{\lambda}(\xi)| \, d\xi$ is finite.  Thus we can define
a function $h(x)$ to be equal to $\int_{{\bf R}^n}
\widehat{\lambda}(\xi) \, \exp (2 \pi i \, x \cdot \xi) \, d\xi$ for
each $x \in {\bf R}^n$, and $h(x)$ is bounded and continuous.
We would like to check that $\lambda$ is defined by integration using
this density $h$.

	Let $f$ be a continuous function on ${\bf R}^n$ with compact
support.  As above $\lambda(f)$ is equal to the limit of
$\lambda(W_\alpha * f)$ as $\alpha \to 0$.  For each $\alpha > 0$,
$\lambda(W_\alpha * f)$ is equal to the integral on ${\bf R}^n$ of
$\lambda * W_\alpha$ times $f$.  We also know that $(\lambda *
W_\alpha)(x)$ can be expressed as the integral of the product of
$\widehat{\lambda}(\xi)$, $G_\alpha(\xi)$, and $\exp(2 \pi i \, x
\cdot \xi)$, where we integrate in $\xi$.  Because we are assuming
that $\widehat{\lambda}$ is integrable, this integral tends to the
integral of $\widehat{\lambda}(\xi)$ times $\exp (2 \pi i \, x \cdot
\xi)$ as $\alpha \to 0$.

	Thus $(\lambda * W_\alpha)(x)$ converges to $h(x)$ as $\alpha
\to 0$, and in fact the convergence is uniform on compact subsets of
${\bf R}^n$.  It follows that the integral of $\lambda * W_\alpha$
times $f$ converges to the integral of $h$ times $f$ as $\alpha \to
0$, and therefore $\lambda(f)$ is equal to the integral of $h$ times
$f$ for all continuous functions $f$ on ${\bf R}^n$ with compact
support.  Using standard arguments one can show that $h$ is integrable
and that $\lambda(f)$ is equal to the integral of $h$ times $f$ for
all bounded continuous functions $f$ on ${\bf R}^n$.

	Now suppose that $h(x)$ is a continuous integrable function on
${\bf R}^n$, which we can view as the density of a finite measure on
${\bf R}^n$.  Thus we can define the Fourier transform of $h$ as the
Fourier transform of that measure, which is to say that
$\widehat{h}(\xi)$ is equal to the integral of $h(x)$ times $\exp(- 2
\pi i \, \xi \cdot x)$ for all $\xi \in {\bf R}^n$.  The convolution
$(h * W_\alpha)(x)$ can be expressed explicitly as the integral of
$h(y)$ times $W_\alpha(x - y)$, where we integrate in $y$, and this
tends to $h(x)$ as $\alpha \to 0$.  It follows that the integral of
$\widehat{h}(\xi)$ times $G_\alpha(\xi)$ times $\exp (2 \pi i \, x
\cdot xi)$, where we integrate in $\xi$, tends to $h(x)$ as $\alpha
\to 0$.  As in the preceding paragraphs, if we assume that
$\widehat{h}(\xi)$ is integrable, then we can simply say that the
integral of $\widehat{h}(\xi)$ times $\exp (2 \pi i \, x \cdot \xi)$
with respect to $\xi$ is equal to $h(x)$ for all $x \in {\bf R}^n$.

	Suppose that $h(x)$ is a continuous integrable function on
${\bf R}^n$ and that the Fourier transform $\widehat{h}(\xi)$ of $h$
is a nonnegative real number for all $\xi \in {\bf R}^n$.  If we put
$x = 0$ in the identities just discussed, we obtain that the integral
of $\widehat{h}(\xi)$ times $G_\alpha(\xi)$ is equal to $(h *
W_\alpha)(0)$ for all $\alpha > 0$.  Because $(h * W_\alpha)(0) \to
h(0)$ as $\alpha \to 0$, we get that the integral of
$\widehat{h}(\xi)$ times $G_\alpha(\xi)$ tends to $h(0)$ as $\alpha
\to 0$.  The hypothesis that $\widehat{h}(\xi)$ is a nonnegative real
number for all $\xi \in {\bf R}^n$ permits us to conclude that
$\widehat{h}$ is integrable on ${\bf R}^n$.

	Let $\lambda$ be a finite measure on ${\bf R}^n$ such that the
Fourier transform $\widehat{\lambda}$ of $\lambda$ has compact support
in ${\bf R}^n$.  Since $\widehat{\lambda}$ is a continuous function,
it follows that $\widehat{\lambda}$ is integrable.  It follows from the
earlier discussion that $\lambda$ corresponds to integration with a
continuous integrable density $h$.

	We also have that $h(x) = \int_{{\bf R}^n}
\widehat{\lambda}(\xi) \, \exp (2 \pi i \, x \cdot \xi) \, d\xi$.
Because $\widehat{\lambda}$ has compact support, the integral of
$\widehat{\lambda}(\xi)$ times $\exp (2 \pi i \, z \cdot \xi)$ with
respect to $\xi$ makes sense for all $z \in {\bf C}^n$.  In fact it
defines a holomorphic function of $z$ on ${\bf C}^n$, which is a
complex analytic extension of $h$ on ${\bf R}^n$.  If $z = x + i y$,
$x, y \in {\bf R}^n$, then at $z$ this extension is bounded by a 
constant times the exponential of a constant times the norm of $y$.

	In particular, if $\widehat{\lambda}$ has compact support,
and if $\lambda$ vanishes on a nonempty open subset of ${\bf R}^n$,
then $\lambda$ is the zero measure.

	Now suppose that $n = 1$, and that $\lambda$ is a finite
measure on the real line such that $\widehat{\lambda}(\xi) = 0$ when
$\xi < 0$.  The integral of $\widehat{\lambda}(\xi)$ times $\exp (2
\pi i \, z \, \xi)$ with respect to $\xi$ makes sense for all complex
numbers $z = x + i y$ with $x, y \in {\bf R}$ and $y > 0$.
This defines a complex analytic function on the upper half plane
in ${\bf C}$.

	In general dimensions, if the Fourier transform of a finite
measure $\lambda$ is supported in some closed set, then one may be
able to make sense of the integral of $\widehat{\lambda}(\xi)$ times
$\exp (2 \pi i \, z \cdot \xi)$ with respect to $\xi$ for some $z \in
{\bf C}^n$, with interesting complex-analyticity properties, and so
on.

\end{document}